\documentclass[11pt,reqno]{amsart}

\theoremstyle{plain}
\newtheorem{theorem} {Theorem}[section]
\newtheorem{lemma}[theorem] {Lemma}

\newtheorem{property}[theorem] {Property}

\theoremstyle{definition}
\newtheorem{definition}[theorem] {Definition}
\newtheorem{conjecture}[theorem] {Conjecture}

\theoremstyle{remark}
\newtheorem{remark}[theorem] {Remark}

\numberwithin{equation}{section}

\usepackage{amsfonts}
\usepackage{amssymb}
\usepackage{amscd}

\newcommand{\R}{{\mathbb R}}
\newcommand{\Z}{{\mathbb Z}}
\newcommand{\N}{{\mathbb N}}

\newcommand{\CC}{{\mathbb C}}

\newcommand{\TT}{{\mathcal T}}

\newcommand{\BB}{{\mathfrak B}}
\newcommand{\al}{{\alpha}}
\newcommand{\la}{{\lambda}}
\newcommand{\sa}{{\sigma}}

\newcommand{\iy}{{\infty}}
\newcommand{\vphi}{{\varphi}}
\newcommand{\vep}{{\varepsilon}}
\newcommand{\g}{{\gamma}}
\newcommand{\de}{{\delta}}

\newcommand{\be}{{\beta}}

\newcommand{\bna}{\begin{eqnarray}}
\newcommand{\ena}{\end{eqnarray}}
\newcommand{\ba}{\begin{eqnarray*}}
\newcommand{\ea}{\end{eqnarray*}}
\newcommand{\beq}{\begin{equation}}
\newcommand{\eeq}{\end{equation}}

\begin{document}

\title[Constants in Bernstein-Nikolskii Inequalities]
{Sharp Constants of Approximation Theory. I. Multivariate Bernstein-Nikolskii Type Inequalities}
\author{Michael I. Ganzburg}
 \address{Department of Mathematics\\ Hampton University\\ Hampton,
 VA 23668\\USA}
 \email{michael.ganzburg@hamptonu.edu}

 \begin{abstract}
 Given a centrally symmetric convex body $V\subset\R^m$, let $\TT_{aV}$ be the set of
 all trigonometric polynomials
 with the spectrum in $aV,\,a>0$, and let $B_V$ be the set of all entire functions of exponential
 type with the spectrum in $V$.
 We discuss limit relations between the sharp constants in the multivariate Bernstein-Nikolskii
inequalities defined by
\ba
&&P_{p,q,D_N,a,V}:=a^{-N-m/p+m/q}\sup_{T\in\TT_{aV}\setminus\{0\}}\frac{\|D_N(T)\|_{L_q\left(Q_{\pi}\right)}}
{\|T\|_{L_p\left(Q_{\pi}\right)}},\\
&&E_{p,q,D_N,V}:=
\sup_{f\in (B_V\cap L_p(\R^m))\setminus\{0\}}\frac{\|D_N(f)\|_{L_q(\R^m)}}
{\|f\|_{L_p(\R^m)}},
\ea
where $0<p\le q\le \iy,\,Q_{\pi}=\{x\in\R^m: \vert x_j\vert\le\pi,\,1\le j\le m\},$
and $ D_N=\sum_{\vert\al\vert=N}b_\al D^\al$ is a differential operator with constant coefficients.
We prove that
\ba
E_{p,q,D_N,V}\le \liminf_{a\to\iy}P_{p,q,D_N,a,V},\qquad
E_{p,\iy,D_N,V}= \lim_{a\to\iy}P_{p,\iy,D_N,a,V}.
\ea
 \end{abstract}
 \maketitle

 \section{Introduction}\label{S1}
\setcounter{equation}{0}
\noindent
In this paper we discuss relations between the sharp constants in the multivariate Bernstein-Nikolskii
type inequalities for trigonometric polynomials and entire functions of exponential type.\vspace{.12in}\\
\textbf{Notation.}
Let $\R^m$ be the Euclidean $m$-dimensional space with elements $x=(x_1,\ldots,x_m),\, y=(y_1,\ldots,y_m),
\,u=(u_1,\ldots,u_m)$, the inner product $x\cdot y:=\sum_{j=1}^mx_jy_j$,
and the norm $\vert x\vert:=\sqrt{x\cdot x}$.
Next, $\CC^m:=\R^m+i\R^m$ is the $m$-dimensional complex space with elements $z=(z_1,\ldots, z_m)$;
$\Z^m$ denotes the set of all integral lattice points in $\R^m$; and $\Z^m_+$ is a subset of $\Z^m$
of all points with nonnegative coordinates.
Given $\sa\in\R^m,\,\sa_j>0,\,1\le j\le m$, and $M>0$, let
$\Pi_\sa:=\{x\in\R^m: \vert x_j\vert\le\sa_j, 1\le j\le m\},\,
Q_M:=\{x\in\R^m: \vert x_j\vert\le M, 1\le j\le m\}$,
and
$\BB_M:=\{x\in\R^m: \vert x\vert\le M\}$
be the $m$-dimensional parallelepiped, cube, and ball, respectively.
In addition, let $d(x,B):=\inf_{y\in B}\vert x-y\vert,\,x\in\R^m$; \,
$d(A,B):=\inf_{x\in A}d(x,B)$; and
 $\mbox{diam}(A):=\sup_{x\in A}\sup_{y\in A}\vert x-y\vert$,
where  $A$ and $B$ are subsets of $\R^m$.

Let $L_r(E)$ be the space of all measurable complex-valued functions $F$
 on a measurable set $E\subseteq\R^m$  with the finite functional
 \ba
 \|F\|_{L_r(E)}:=\left\{\begin{array}{ll}
 \left(\int_E\vert F(x)\vert^r dx\right)^{1/r}, & 0<r<\iy,\\
 \mbox{ess} \sup_{x\in E} \vert F(x)\vert, &r=\iy.
 \end{array}\right.
 \ea
 This functional allows the following "triangle" inequality:
 \beq\label{E1.1a}
 \left\|\sum_{j=1}^s F_j\right\|^{\tilde{r}}_{L_r(E)}
 \le \sum_{j=1}^s \left\|F_j\right\|^{\tilde{r}}_{L_r(E)},
 \qquad F_j\in L_r(E),\qquad
 1\le j\le s,
 \eeq
 where $s\in\N:=\{1,\,2,\ldots\}$ and $\tilde{r}:=\min\{1,r\}$ for $r\in(0,\iy]$.

 Throughout the paper $V$ is a centrally symmetric (with respect to the origin)
 closed
 convex body in $\R^m$. Its $m$-dimensional volume is denoted by $\vert V\vert_m$.
 The set $V$ generates the following dual norms
 on $\R^m$ and $\CC^m$ by
 \ba
 \|y\|_V^*:=\sup_{x\in V}\vert x\cdot y\vert,\quad y\in\R^m;\qquad
 \|z\|_V^*:=\sup_{x\in V}\left\vert\sum_{j=1}^m x_jz_j\right\vert,\quad z\in\CC^m.
 \ea
 For example, if $\sa\in\R^m,\,\sa_j>0,\,1\le j\le m$, and $V=\left\{x\in\R^m:
 \left(\sum_{j=1}^m\vert x_j/\sa_j\vert^{\mu}\right)^{1/\mu}\le 1\right\}$, then
 for $y\in\R^m$,
 $\|y\|_V^*=\left(\sum_{j=1}^m\vert \sa_j y_j\vert^\la\right)^{1/\la}$,
 where $\mu\in[1,\iy],\,\la\in[1,\iy]$, and $1/\mu+1/\la=1$.
 In particular,
 $\|y\|_{\Pi_\sa}^*=\sum_{j=1}^m\sa_j\vert y_j\vert,\,
 \|y\|_{Q_M}^*=M\sum_{j=1}^m\vert y_j\vert$,
 and $\|y\|_{\BB_M}^*=M\vert y\vert$.

 Given $a>0$, the set of all trigonometric polynomials
 $T(x)=\sum_{k\in aV\cap \Z^m}c_ke^{ik\cdot x}$
 with complex coefficients is denoted by $\TT_{aV}$.

  \begin{definition}\label{D1.1}
 We say that an entire function $f:\CC^m\to \CC^1$ has exponential type $V$
 if for any $\vep>0$ there exists a constant $C_0(\vep,f)$ such that
 for all $z\in\CC^m,\,
 \vert f(z)\vert\le C_0(\vep,f)\exp\left((1+\vep)\|z\|_V^*\right)$.
 \end{definition}
  The class of all entire function of exponential type $V$ is denoted
  by $B_V$.
  Throughout the paper, if no confusion may occur, the same notation is applied to
  $f\in B_V$ and its restriction to $\R^m$ (e. g., in the form
  $f\in  B_V\cap L_p(\R^m)$).
  The class $B_V$ was defined by Stein and Weiss
  \cite[Sect. 3.4]{SW1971}. For $V=\Pi_\sa,\,V=Q_M,$ and $V=\BB_M$, similar
  classes were
  defined by Bernstein \cite{B1948} and Nikolskii
  \cite{N1951}, \cite[Sects. 3.1, 3.2.6]{N1969}, see also
  \cite[Definition 5.1]{DP2010}.
  Properties of functions from $B_V$ have been investigated in numerous
  publications (see, e. g., \cite{B1948, N1951, N1969, SW1971, NW1978, G2001} and
  references therein). Some of these properties are presented in Section \ref{S2a}.
  In particular, if $f\in  B_V\cap L_\iy(\R^m)$, then the norm $\|z\|_V^*$
  in Definition \ref{D1.1} can be replaced with
  $\|\mbox{Im}\,z\|_V^*$ (see Lemma \ref{L2.1}(b)).

   The Fourier transform of a function $\psi\in L_1(\R^m)$ or $\psi\in L_2(\R^m)$
   is denoted by the formula
   \ba
   \widehat{\psi}(y):=
   \frac{1}{(2\pi)^{m/2}}\int_{\R^m}\psi(x)e^{-ix\cdot y}dx.
   \ea
\noindent
We use the same notation $\widehat{\psi}$
for the Fourier transform of a tempered distribution
 $\psi$ on $\R^m$.
 By the definition (see, e. g., \cite[Sect. 1.3]{SW1971}),
 $\psi$ is a continuous linear functional $\langle \psi,\vphi\rangle$ on the Schwartz
 class $S(\R^m)$ of all
 test functions $\vphi$ on $\R^m$, and $\widehat{\psi}$ is defined by the formula
 $\langle \widehat{\psi},\vphi\rangle
 :=\langle\psi,\widehat{\vphi}\rangle,\,\vphi\in S(\R^m).$
 The support of a function or a distribution $\psi$ on $\R^m$ is denoted by
 $\mbox{supp}(\psi)$.

Throughout the paper $C,\,C_0,\,C_1,\ldots$ denote positive constants independent
of essential parameters.
 Occasionally we indicate dependence on certain parameters.

 We also use multi-indices $\al=(\al_1,\ldots,\al_m)\in \Z^m_+$ and
 $\be=(\be_1,\ldots,\be_m)\in \Z^m_+$ with
 $\vert\al\vert:=\sum_{j=1}^m\al_j,\,\be\le \al$
 (i. e., $\be_j\le\al_j,\,1\le j\le m$), and
 \ba
 \binom{\al}{\be}:=\frac{\prod_{j=1}^m\al_j!}
 {\prod_{j=1}^m\be_j!\prod_{j=1}^m(\al_j-\be_j)!},\quad
 x^\al:=x_1^{\al_1}\cdot\cdot\cdot x_m^{\al_m}, \quad
 D^\al:=\frac{\partial^{\al_1}}{\partial x_1^{\al_1}}\cdot\cdot\cdot
 \frac{\partial^{\al_m}}{\partial x_m^{\al_m}}.
 \ea
 \noindent
 In addition, we use the ceiling and floor functions
 $\lceil a \rceil$ and $\lfloor a \rfloor$.
  \vspace{.12in}\\
\textbf{Multivariate Bernstein-Nikolskii Type Inequalities.}
Let
$D_N:=\sum_{\vert\al\vert=N}b_\al D^\al$
be a linear differential
operator with constant coefficients $b_\al\in\CC^1,\,\vert\al\vert=N$,
and let $\Delta_N(y):=\sum_{\vert\al\vert=N}b_\al y^\al,\,y\in\R^m,$
be the total symbol of $D_N,\,
N\in \N$.
We assume that $D_0$
is the identity operator and $\Delta_0:=1$.

 Next, we define  sharp constants in multivariate Bernstein-Nikolskii type
inequalities for trigonometric polynomials and entire functions of exponential type. Let
\bna
&&P_{p,q,D_N,a,V}:=a^{-N-m/p+m/q}
\sup_{T\in\TT_{aV}\setminus\{0\}}\frac{\|D_N(T)\|_{L_q(Q_\pi)}}
{\|T\|_{L_p(Q_\pi)}},\label{E1.1n}\\
&& E_{p,q,D_N,V}:=
\sup_{f\in (B_{V}\cap L_p(\R^m))\setminus\{0\}}\frac{\|D_N(f)\|_{L_q(\R^m)}}
{\|f\|_{L_p(\R^m)}}.\label{E1.2n}
\ena
Here, $a>0,\,N\in\Z^1_+,\,V\subset\R^m$, and $0<p\le q\le\iy$.
Note that
\beq\label{E1.3n}
E_{p,q,D_N,\g V}=\g^{N+m/p-m/q}\,E_{p,q,D_N,V},\qquad \g>0.
\eeq
Indeed, by Definition \ref{D1.1}, $g\in B_{\g V}$ if and only if $f(x)=g(x/\g)\in B_V$,
and \eqref{E1.3n} follows from the relations
\ba
E_{p,q,D_N,\g V}
&=&\g^{N+m/p-m/q}
\sup_{g\in (B_{\g V}\cap L_p(\R^m))\setminus\{0\}}\frac{\|D_N(g(\g^{-1}\cdot))\|_{L_q(\R^m)}}
{\|g(\g^{-1}\cdot)\|_{L_p(\R^m)}}\\
&=&\g^{N+m/p-m/q}
\sup_{f\in (B_{V}\cap L_p(\R^m))\setminus\{0\}}\frac{\|D_N(f)\|_{L_q(\R^m)}}
{\|f\|_{L_p(\R^m)}}.
\ea

A detailed survey of univariate Bernstein-Nikolskii  inequalities for
$V=[-\sa,\sa],\,D_N=d^N/dx^N,\,a=n/\sa,\,\sa>0,\,n\in\N$,
and $0<p\le q\le\iy$ was presented in \cite{GT2017},
so we assume that $m>1$ in all discussions below.

The Bernstein-Nikolskii sharp constants $P_{p,q,D_N,a,V}$ and $E_{p,q,D_N,V}$ for $0<p< q\le\iy$
can be easily found only for $p=2$ and $q=\iy$. Namely,
\bna
&&P_{2,\iy,D_N,a,V}
=(2\pi)^{-m/2} a^{-(N+m/2)}
\left(\sum_{k\in aV\cap\Z^m}\vert \Delta_N(i k)\vert^2\right)^{1/2},
\label{E1.4n}\\
&&E_{2,\iy,D_N,V}=(2\pi)^{-m/2} \left(\int_V\vert \Delta_N(i x)\vert^2\,dx\right)^{1/2}.
\label{E1.5n}
\ena
For $N=0$ these equalities were obtained by Nessel and Wilmes \cite[pp. 10, 11]{NW1978}.
The Bernstein sharp constants for $0<p= q\le\iy$ can be found in some special cases as well.
In particular,
\bna
&&P_{2,2,D_N,a,V}=a^{-N}\max_{k\in aV\cap\Z^m}\vert \Delta_N(i k)\vert,\quad
E_{2,2,D_N,V}=\max_{x\in V}\vert \Delta_N(i x)\vert, \label{E1.6n}\\
&&P_{q,q,D^\al,a,\Pi_\sa}=\sa^\al(1+o(1)),\quad
E_{q,q,D^\al,\Pi_\sa}=\sa^\al,\quad \vert\al\vert=N,\quad q\in(0,\iy],\label{E1.7n}\\
&&P_{q,q,D^\al,a,Q_M}=M^{\vert\al\vert}(1+o(1)),\quad
E_{q,q,D^\al,Q_M}=M^{\vert\al\vert},\quad \vert\al\vert=N,\quad q\in(0,\iy],\label{E1.7na}
\ena
as $a\to \iy$.
Note that the following limit relations hold true:
\bna
&&E_{2,\iy,D_N,V}=\lim_{a\to\iy}P_{2,\iy,D_N,a,V},\label{E1.7an}\\
&&E_{2,2,D_N,V}=\lim_{a\to\iy}P_{2,2,D_N,a,V},\label{E1.7bn}\\
&&E_{q,q,D^\al,\Pi_\sa}=\lim_{a\to\iy}P_{q,q,D^\al,a,\Pi_\sa},\label{E1.7cn}\\
&&E_{q,q,D^\al,Q_M}=\lim_{a\to\iy}P_{q,q,D^\al,a,Q_M}.\label{E1.7dn}
\ena
Simple proofs of relations \eqref{E1.4n}-\eqref{E1.7bn} are presented in Lemmas
\ref{L2.4}(a), \ref{L2.5}, and \ref{L2.6},
while \eqref{E1.7cn} and \eqref{E1.7dn} follow immediately from \eqref{E1.7n}
and \eqref{E1.7na}, respectively.
Certain sharp Bernstein-type inequalities for functions from $B_V$ in case of
$p=q=\iy$ are known for special $V$ and $D_N$.
In particular, Kamzolov \cite[Corollary 2]{K1974} proved the relation
\beq\label{E1.7nn}
E_{\iy,\iy,\Delta,\BB_M}=m M^2,
\eeq
where $\Delta$ is the Laplace operator.
The author \cite[Theorem 5]{G1979},
\cite[Theorem 2]{G1982} extended \eqref{E1.7nn} to
elliptic differential operators $D_2$ with constant
coefficients and $m$-dimensional ellipsoids $V$.

So there are just a few examples of finding exact or asymptotically
exact values of $P_{p,q,D_N,a,V}$ and $E_{p,q,D_N,V}$ for $0<p\le q\le\iy$.
However, efficient estimates of these constants are possible.

The following multivariate Nikolskii-type inequalities for $0<p\le q\le\iy$
and functions from
$\TT_{aV}$ and $B_V$ were proved by Nessel and Wilmes
\cite[Theorems 2, 5]{NW1978}:
\bna
&&P_{p,q,D_0,a,V}\le \left[(\lceil p/2\rceil/(2\pi))^m\vert V\vert_m
\right]^{1/p-1/q}+o(1),
\qquad a\to\iy,\label{E1.8n}\\
&&E_{p,q,D_0,V}\le \left[(\lceil p/2\rceil/(2\pi))^m\vert V\vert_m
\right]^{1/p-1/q}.
\label{E1.9n}
\ena
For $V=\Pi_\sa$ and $1\le p\le q\le\iy$ similar inequalities were established
 by Nikolskii \cite{N1951}, \cite[Sects. 3.3.5, 3.4.3]{N1969}.

Combining estimates \eqref{E1.8n} and \eqref{E1.9n} with \eqref{E1.7na},
we arrive at the crude Bernstein-Nikolskii type inequalities
($a>0,\,N\in\Z^1_+,\,V\subset\R^m,\,0<p\le q\le\iy$)
\bna
&&P_{p,q,D_N,a,V}\le C_1 (\mbox{diam}(V))^{N+m/p-m/q}+o(1),\qquad a\to\iy,\label{E1.10n}\\
&&E_{p,q,D_N,V}\le C_2 (\mbox{diam}(V))^{N+m/p-m/q},\label{E1.11n}
\ena
where $C_1$ is independent of $a$ and $V$, and $C_2$ is independent of $V$
 (see Lemma \ref{L2.4}(b) for a short proof).
\vspace{.12in}\\
\textbf{Main Results.}
 Our major results discuss relations between $P_{p,q,D_N,a,V}$ and
 $E_{p,q,D_N,V}$.
 In particular, we extend
 relation \eqref{E1.7an} to any $p\in(0,\iy)$.

 \begin{theorem} \label{T1.2}
 If $a>0,\,N\in\Z^1_+,\,V\subset\R^m$, and $0<p\le q\le\iy$,
 then the following relation holds true:
 \beq
 E_{p,q,D_N,V}\le \liminf_{a\to\iy}P_{p,q,D_N,a,V}.\label{E1.12n}\\
\eeq
 \end{theorem}
 \noindent
 In case of $q=\iy$ a more precise result is valid.
 \begin{theorem} \label{T1.3}
 If $a>0,\,N\in\Z^1_+,\,V\subset\R^m$, and $p\in(0,\iy]$, then
  $ \lim_{a\to\iy}P_{p,\iy,D_N,a,V}$ exists and
 \beq \label{E1.13n}
  E_{p,\iy,D_N,V}=\lim_{a\to\iy}P_{p,\iy,D_N,a,V}.
 \eeq
 In addition, there exists a nontrivial function $f_0\in  B_V\cap L_p(\R^m)$ such that
 \beq \label{E1.14n}
 \|D_N(f_0)\|_{L_\iy(\R^m)}/\|f_0\|_{L_p(\R^m)}=\lim_{a\to\iy}P_{p,\iy,D_N,a,V}.
 \eeq
 \end{theorem}

 \begin{remark}\label{R1.4}
Relations \eqref{E1.13n} and \eqref{E1.14n} show that the function
$f_0\in B_V\cap L_p(\R^m)$
from Theorem \ref{T1.3} is an extremal function for $E_{p,\iy,D_N,V}$.
 \end{remark}

 \begin{remark}\label{R1.5}
 In view of \eqref{E1.13n}, \eqref{E1.7bn}, and \eqref{E1.7cn}
 we believe that the following conjecture is valid.
 \begin{conjecture}\label{C1.6}
 The limit $\lim_{a\to\iy}P_{p,q,D_N,a,V}$ exists and
 $E_{p,q,D_N,V}=\lim_{a\to\iy}P_{p,q,D_N,a,V}$  for
 $0<p\le q<\iy$.
 \end{conjecture}
 \end{remark}

\begin{remark}\label{R1.6}
In definitions \eqref{E1.1n} and \eqref{E1.2n} of the sharp constants we
discuss only complex-valued functions $T$ and $f$. We can define similarly
the "real" sharp constants if the suprema in \eqref{E1.1n} and \eqref{E1.2n}
 are taken over all real-valued functions
on $\R^m$ from $\TT_{aV}\setminus\{0\}$ and $(B_V\cap L_p(\R^m))\setminus\{0\}$.

We do not know as to whether the "complex" and "real" sharp constants coincide
but for $q=\iy$ this is true. For $m=1$ this fact was proved in
\cite[Theorem 1.1]{GT2017} and the case of $m>1$ can be proved similarly.

"Real" analogues of Theorems \ref{T1.2} and \ref{T1.3} are valid as well.
The proof of the "real" version of \eqref{E1.12n} is similar to that of
Theorem \ref{T1.2} if we use Property \ref{P2.1} from Section \ref{S2}.
The "real" version of \eqref{E1.13n} follows immediately from the fact
mentioned above.
The "real" version of Conjecture \ref{C1.6} is believed to be true as well.
\end{remark}

\begin{remark}\label{R1.7}
In the univariate case of $V=[-1,1],\,D_N=d^N/dx^N$, and $a\in\N$,
 Theorems \ref {T1.2} and \ref{T1.3} for the "real" and
"complex" sharp constants were proved by the author and Tikhonov \cite{GT2017}.
In earlier publications \cite{LL2015a, LL2015b}, Levin and Lubinsky established
versions of Theorems \ref {T1.2} and \ref{T1.3} on the unit circle for $N=0$.
More precise asymptotics for $P_{p,\iy,D_0,a,[-1,1]}$ were obtained
by Gorbachev and Martyanov \cite[Theorem 1]{GM2018}.
Certain extensions of the Levin-Lubinsky's results to the $m$-dimensional
unit sphere in $\R^{m+1}$ were recently proved by Dai, Gorbachev, and Tikhonov
\cite[Theorem 1.1]{DGT2017}.
\end{remark}

The proofs of Theorems \ref {T1.2} and \ref{T1.3} are presented in Section \ref{S4}.
Section \ref{S2a} contains certain properties of functions from $B_V$ and $\TT_{aV}$.
Multivariate Levitan's polynomials are introduced in Section \ref{S2}.

\section{Properties of Entire Functions and Trigonometric
  Polynomials}\label{S2a}
 \noindent
\setcounter{equation}{0}
\noindent
In this section we discuss certain properties of functions from $B_V$ and $\TT_{aV}$
that are needed for the proofs of Theorems \ref {T1.2} and \ref{T1.3} (see
Lemmas \ref{L2.1}, \ref{L2.2}, \ref{L2.2a}).
In addition, we prove
here certain multivariate Bernstein-Nikolskii type inequalities
presented in Section \ref {S1} (see
Lemmas \ref{L2.3} through \ref{L2.6}).
Certain standard facts are included in the following lemma.

\begin{lemma}\label{L2.1}
(a) If $V_1\subseteq V_2$, then $\TT_{aV_1}\subseteq\TT_{aV_2}$
and  $B_{V_1}\subseteq B_{V_2}$.\\
(b) If $f\in B_V\cap L_\iy(\R^m)$, then for every $x\in\R^m$ and $y\in\R^m$,
\ba
\vert f(x+iy)\vert\le \|f\|_{L_\iy(\R^m)}\exp\left(\|y\|_V^*\right).
\ea
(c) Let $f$ be a tempered distribution on $\R^m$.
\begin{itemize}
\item[(i)] If $f\in B_V$, then $\mathrm{supp}\left(\widehat{f}\right)\subseteq V$.
\item[(ii)] If $\mathrm{supp}\left(\widehat{f}\right)\subseteq V$, then $f$ can be
extended to $\CC^m$ as a function from $B_V$.
\end{itemize}
(d) Let $f$ be a tempered distribution on $\R^m$.
If $f\in B_V$, then $D^\al f\in B_V$ for every $\al\in \Z^m_+$.
\end{lemma}
\proof
Statement (a) of the lemma follows immediately from the definitions
of $\TT_{aV}$ and $B_V$.
Statement (b) is established in \cite[Eq. (4.13)]{NW1978}
(cf. \cite[Lemma 3.4.11]{SW1971}), while the proof of a Paley-Wiener-Schwartz
 type theorem (c) is outlined
in \cite[p. 13]{NW1978}.
It remains to prove statement (d). Let $\vphi$ be a test function from $S(\R^m)$
with $\mbox{supp}(\vphi) \subseteq\R^m\setminus V$.
Then for any $\al\in\Z^m_+$, the function $\vphi_1(x):=x^\al\vphi(x)$ belongs to
$S(\R^m)$ and $\mbox{supp}(\vphi_1) \subseteq\R^m\setminus V$.
Therefore, by part (i) of statement (c),
\ba
\left\langle\widehat{D^\al f},\vphi\right\rangle=i^{\vert\al\vert}
\left\langle\widehat{ f},\vphi_1\right\rangle=0.
\ea
So $\mbox{supp}(\widehat{D^{\al}f}) \subseteq V$, and by part (ii) of statement (c),
$D^\al f\in B_V$.\hfill $\Box$

\begin{remark}\label{R2.1a}
The condition in statements (c) and (d) of Lemma \ref{L2.1}
that $f$ is a tempered distribution on $\R^m$ is obviously satisfied for
$f\in L_p(\R^m),\,p\in(0,\iy]$. In particular, when $f\in L_2(\R^m)$,
a Paley-Wiener type theorem of Lemma \ref{L2.1} (c) was proved
in \cite[Theorem 3.4.9]{SW1971}.
\end{remark}
\noindent
The compactness theorem for the set $B_V\cap L_\iy(\R^m)$ is discussed below.

\begin{lemma}\label{L2.2}
For any sequence $\{f_n\}_{n=1}^\iy,\,
f_n\in B_V\cap L_\iy(\R^m),\,n\in\N,$
with $\sup_{n\in\N}\| f_n\|_{L_\iy(\R^m)}= C$, there exist a subsequence
$\{f_{n_s}\}_{s=1}^\iy$ and a function $f_0\in B_V\cap L_\iy(\R^m)$
such that for every $\al\in\Z^m_+$,
\beq\label{E2.1}
\lim_{s\to\iy} D^\al f_{n_s}=D^\al f_0
\eeq
uniformly on any compact set in $\CC^m$.
\end{lemma}
\proof
For $V=\Pi_\sa$ this compactness result was proved by Nikolskii
\cite[Theorem 3.3.6]{N1969}. Next, for any $V\subset \R^m$
there exists $\sa(V)\in\R^m,\,\sa(V)_j>0,\,1\le j\le m$, such that
$V\subseteq \Pi_{\sa(V)}$.
Note that by Lemma \ref{L2.1}(a),
$f_n\in B_{\Pi_{\sa(V)}}\cap L_\iy(\R^m),\,n\in\N$.
Then by the Nikolskii's compactness theorem,
there exists a subsequence
$\{f_{n_s}\}_{s=1}^\iy$ and a function $f_0\in B_{\Pi_{\sa(V)}}\cap L_\iy(\R^m)$
such that \eqref{E2.1} holds true uniformly on any compact set in $\CC^m$.
It remains to show that $f_0\in B_V\cap L_\iy(\R^m)$.
Indeed, by Lemma \ref{L2.1}(b),
\beq\label{E2.2}
\vert f_n(x+iy)\vert\le C\exp\left(\|y\|_V^*\right),\qquad n\in\N,\quad x\in\R^m,
\quad y\in\R^m.
\eeq
Therefore using first \eqref{E2.1} for $\al=0$ and then \eqref{E2.2}, we obtain
for any $x\in\R^m$ and $y\in\R^m$
\ba
\vert f_0(x+iy)\vert
=\lim_{s\to\iy}\vert f_{n_s}(x+iy)\vert
\le C\exp\left(\|y\|_V^*\right)
\le C\exp\left(\|x+iy\|_V^*\right).
\ea
Thus $f_0\in B_V\cap L_\iy(\R^m)$. \hfill $\Box$\vspace{.12in}\\
\noindent
Periodization properties of functions from $B_V$ are based on
the following version of the Poisson summation formula
(cf. \cite[Sect. 7.2]{SW1971}).

\begin{lemma}\label{L2.2a}
If $g\in L_1(\R^m)$ and the series
$\sum_{k\in\Z^m}\widehat{g}(k) e^{ik\cdot y}$
 is a finite sum, then the series
 $\sum_{k\in\Z^m}g(y+2k\pi)$ converges in $L_1(Q_\pi)$
 to a function from $L_1(Q_\pi)$ and
 \beq\label{E2.2a}
 \sum_{k\in\Z^m}g(y+2k\pi)
 =(2\pi)^{-m/2}\sum_{k\in\Z^m}\widehat{g}(k) e^{ik\cdot y}.
 \eeq
\end{lemma}
\proof
It was proved in \cite[Theorem 7.2.4]{SW1971}) that
if  $g\in L_1(\R^m)$, then
 $\sum_{k\in\Z^m}g(y+2k\pi)$ converges in $L_1(Q_\pi)$
 to a function $\tilde{g}\in L_1(Q_\pi)$ and its Fourier
 series expansion coincides with the right-hand side
 of \eqref{E2.2a}.
 Since it is a trigonometric polynomial $T$, we
 conclude that $\tilde{g}=T$.\hfill $\Box$

\begin{lemma}\label{L2.3}
For any $M>0$ and $q\in(0,\iy]$ there exists a family
$\{f_{h,q,M}\}_{h\in(0,M)}=\{f_h\}_{h\in(0,M)}$ of
functions from $B_{[-M,M]}$ such that for any $s\in\N$,
\ba
\frac{\left\|f_h^{(s)}\right\|_{L_q(\R^1)}}
{\left\|f_h\right\|_{L_q(\R^1)}}\ge M^s(1-o(1)),\qquad h\to 0^+.
\ea
\end{lemma}
\proof
For $q\in[1,\iy]$ and $s=1$ this result was proved by Akhiezer
\cite[Sect. 84]{A1965}. We will use the similar construction.
Let $d:=\lfloor 1/q\rfloor +1$ and let $\vphi$ be a continuously
$d$-differentiable function on $[0,1]$, satisfying the following
boundary conditions: $\vphi^{(l)}(0)=\vphi^{(l)}(1)=0,\,0\le l\le d$.
Then the function
\ba
f_h(t):=e^{iMt}\int_0^1e^{-iht\tau}\vphi(\tau)\,d\tau
=\frac{1}{h}\int_{M-h}^{M}
e^{it\tau}\vphi\left(\frac{M-\tau}{h}\right)\,d\tau,
\quad h\in(0,M),
\ea
belongs to $B_{[-M,M]}\cap L_q(\R^1)$ because integration
 by parts in the first integral shows that
 \ba
 \vert f_h(t)\vert\le \min\left\{\int_0^1\vert\vphi(\tau)\vert\,d\tau,
 \frac{1}{(ht)^d}\int_0^1\left\vert\vphi^{(d)}(\tau)\right\vert\,d\tau\right\},
 \quad t\in\R^1.
 \ea
 Next,
 \ba
 f_h^{(s)}(t)=i^se^{iMt}\sum_{l=0}^s(-1)^{s-l}\binom{s}{l}M^l
 h^{s-l}\psi_{s-l}(ht),
 \ea
 where
 \ba
 \psi_r(t):=\int_0^1e^{-it\tau}\tau^r\vphi(\tau)\,d\tau,\qquad 0\le r\le s.
 \ea
 Note that $\psi_r\in B_{[-1,1]}\cap L_q(\R^1),\,0\le r\le s$, by
 Bernstein's inequality, since
 $\vert \psi_0(ht)\vert=\vert f_h(t)\vert$ and
 $\vert \psi_r(t)\vert=\left\vert \psi_0^{(r)}(t)\right\vert$.
 Then we have by \eqref{E1.1a}
 \ba
 \frac{\left\|f^{(s)}_h\right\|^{\tilde{q}}_{L_q(\R^1)}}
 {\|f_h\|^{\tilde{q}}_{L_q(\R^1)}}
 \ge M^{s\tilde{q}}-h^{\tilde{q}}\sum_{l=0}^{s-1}\binom{s}{l}^{\tilde{q}}
 M^{l\tilde{q}}
 h^{(s-l-1)\tilde{q}}
 \frac{\left\|\psi_{s-l}\right\|^{\tilde{q}}_{L_q(\R^1)}}
 {\|\psi_0\|^{\tilde{q}}_{L_q(\R^1)}}
 = M^{s\tilde{q}}(1-o(1))^{\tilde{q}},\quad h\to 0^+.
\ea
This establishes the lemma.\hfill $\Box$

\begin{lemma}\label{L2.4}
(a) Relations \eqref{E1.7n} and \eqref{E1.7na} hold true.
(b) Inequalities \eqref{E1.10n} and \eqref{E1.11n} are valid.
\end{lemma}
\proof
(a) We first note that in the univariate case
\beq\label{E2.3}
E_{q,q,d^{s}/dx^s,[-M,M]}=M^s,\qquad q\in(0,\iy],\quad s\in\N.
\eeq
Indeed, the proof of the classical inequality
$E_{q,q,d^{s}/dx^s,[-M,M]}\le M^s$ can be found in
\cite[Theorem 11.3.3]{B1954} for $q\in[1,\iy]$ and in \cite{RS1990}
for $q\in(0,1)$.
Then \eqref{E2.3} follows from Lemma \ref{L2.3}. Next, it follows from
\eqref{E2.3} that for every $f\in B_{\Pi_\sa}\cap L_q(\R^m),\,q\in(0,\iy]
,\,\al_j\in\Z^1_+$, and
$1\le j\le m$,
\beq\label{E2.4}
\left(\int_{R^1}\left\vert \frac{\partial^{\al_j}f(x_1,\ldots,x_j,\ldots,x_m)}
{\partial x_j^{\al_j}}\right\vert ^qdx_j\right)^{1/q}
\le\sa_j^{\al_j}
\left(\int_{\R^1}\vert f(x_1,\ldots,x_j,\ldots,x_m)\vert^qdx_j\right)^{1/q}.
\eeq
Therefore,
\beq\label{E2.5}
\left\|\frac{\partial^{\al_j}f}
{\partial x_j^{\al_j}}\right\|_{L_q(\R^m)}
\le \sa_j^{\al_j}\|f\|_{L_q(\R^m)},\qquad 1\le j\le m.
\eeq
Using this inequality $m$ times, we arrive at the inequality
\beq\label{E2.6}
E_{q,q,D^\al,\Pi_\sa}\le \sa^\al.
\eeq
This inequality is well-known for $q\in[1,\iy]$
 (see, e. g.,
\cite[Eq. (3.2.2.8)]{N1969}).
Finally, let $f_{h,q,M}$ be functions from Lemma \ref{L2.3}.
Then
the functions
$F_{h}:=\prod_{j=1}^mf_{h,q,\sa_j}(x_j),\,h\in(0,\min_{1\le j\le m}\sa_j)$,
satisfy the inequality
 \beq\label{E2.7}
 \left\|D^\al F_h\right\|_{L_q(\R^m)}
 \ge \sa^\al (1-o(1))\left\|F_h\right\|_{L_q(\R^m)},
 \qquad h\to 0^+,\quad q\in(0,\iy].
 \eeq
  Thus
the second equality in \eqref{E1.7n} follows from
\eqref{E2.6} and \eqref{E2.7}.

It is well known that a periodic analog of \eqref{E2.3}
for $aM\in\N$ is
\beq\label{E2.8}
P_{q,q,d^{s}/dx^s,a,[-M,M]}=M^s,\qquad q\in(0,\iy],\quad s\in\N,
\eeq
see, e. g., \cite[Sect. 4.8.62]{T1963} for $q\in[1,\iy]$
 and see \cite{A1981}
for $q\in(0,1)$. An extremal polynomial in \eqref{E2.8} is
$\cos(aM\tau)$. Then similarly to
\eqref{E2.3} - \eqref{E2.6}, we obtain from \eqref{E2.8}
\beq\label{E2.9}
P_{q,q,D^\al,a,\Pi_\sa}
\le a^{-\vert \al \vert}
\prod_{j=1}^m\lceil a\sa_j\rceil^{\al_j}=\sa^\al(1+o(1))
,\qquad a\to\iy.
\eeq
Finally, the polynomial
$T(x):=\prod_{j=1}^m\cos\left(\lfloor a\sa_j\rfloor x_j\right)$
satisfies the inequality
\beq\label{E2.10}
P_{q,q,D^\al,a,\Pi_\sa}
\ge a^{-\vert \al \vert}
\frac{\left\|D^\al T\right\|_{L_q(Q_\pi)}}{\left\|T\right\|_{L_q(Q_\pi)}}
=a^{-\vert \al \vert}
\prod_{j=1}^m\left(\lfloor a\sa_j\rfloor\right)^{\al_j}
=\sa^\al(1+o(1)),\quad a\to\iy.
\eeq
Thus the first equality in \eqref{E1.7n} follows from
\eqref{E2.9} and \eqref{E2.10}. Relations \eqref{E1.7na} follow
immediately from \eqref{E1.7n}.  This completes the proof of statement (a).
\vspace{.12in}\\
(b) Setting $M=\mbox{diam}(V)/2$, we see
that $V\subseteq \BB_M\subseteq Q_M$.
In particular, $\vert V\vert_m\le \vert \BB_1\vert_m (\mbox{diam}(V)/2)^m$.
 Then using \eqref{E1.1a}, Lemma \ref{L2.1}(a), and the first relation
 of \eqref{E1.7na} (see Lemma \ref{L2.4}(a) for the proof), we obtain
\bna\label{E2.11}
P_{q,q,D_N,a,V}
&\le& P_{q,q,D_N,a,Q_M}
\le \left(\sum_{\vert\al\vert=N}\vert b_\al\vert^{\tilde{q}}
P_{q,q,D^\al,a,Q_M}^{\tilde{q}}\right)^{1/\tilde{q}}\nonumber\\
&=& \left(\frac{\mbox{diam}(V)}{2}\right)^N\left(\sum_{\vert\al\vert=N}\vert b_\al\vert^{\tilde{q}}
\right)^{1/\tilde{q}} +o(1),
\ena
as $a\to\iy$. Taking account of the relations
\ba
P_{p,q,D_N,a,V}
&=&
\sup_{T\in\TT_{aV}\setminus\{0\}}\left(a^{-N}\frac{\|D_N(T)\|_{L_q(Q_\pi)}}
{\|T\|_{L_q(Q_\pi)}}\right)
\left(a^{-m/p+m/q}\frac{\|T\|_{L_q(Q_\pi)}}
{\|T\|_{L_p(Q_\pi)}}\right)\\
&\le & P_{q,q,D_N,a,V} P_{p,q,D_0,a,V},
\ea
we arrive at \eqref{E1.10n} from \eqref{E2.11} and \eqref{E1.8n}.
Inequality \eqref{E1.11n} can be proved similarly.\hfill $\Box$

\begin{lemma}\label{L2.5}
Relations \eqref{E1.4n}, \eqref{E1.5n}, and \eqref{E1.6n} hold true.
\end{lemma}
\proof
The proofs of \eqref{E1.5n} and the second relation in \eqref{E1.6n}
are based on the Paley-Wiener type theorem (see Lemma \ref{L2.1}(c)
and Remark \ref{R2.1a}).
Let $f\in B_V\cap L_2(\R^m)$. Then
\ba
f(u)=\frac{1}{(2\pi)^{m/2}}\int_{V}g(x)e^{iu\cdot x}dx,\quad
D_N(f)(u)=\frac{1}{(2\pi)^{m/2}}\int_{V}g(x)\Delta_N(ix)e^{iu\cdot x}dx,
\ea
where $\|f\|_{L_2(\R^m)}=\|g\|_{L_2(V)}$
and $\|D_N(f)\|_{L_2(\R^m)}=\|g \Delta_N(i\cdot)\|_{L_2(V)}$.
Therefore,
\ba
E_{2,2,D_N,V}
=\sup_{g\in L_2(V)\setminus\{0\}}
\frac{\left\|g\,\Delta_N(i\cdot)\right\|_{L_2(V)}}
{\|g\,|_{L_2(V)}}
=\left\|\Delta_N(i\cdot)\right\|_{L_\iy(V)}.
\ea
Next, it is easy to see that
\ba
E_{2,\iy,D_N,V}
&=&\sup_{f\in (B_V\cap L_2(V))\setminus\{0\}}
\frac{\left\vert D_N(f)(0)\right\vert}
{\|f\,\|_{L_2(\R^m)}}\\
&=&\frac{1}{(2\pi)^{m/2}}\sup_{g\in L_2(V)\setminus\{0\}}
\frac{\left\vert \int_V g(x)\Delta_N(ix)\,dx\right\vert}
{\|g\,\|_{L_2(V)}}\\
&=&\frac{1}{(2\pi)^{m/2}}\left\|\Delta_N(i\cdot)
\right\|_{L_2(V)}.
\ea
The corresponding results for trigonometric polynomials can be
 proved similarly by using Parseval's identity.\hfill $\Box$

\begin{lemma}\label{L2.6}
Limit relations \eqref{E1.7an} and \eqref{E1.7bn} hold true.
\end{lemma}
\proof
Let $\left\|\Delta_N(i\cdot)\right\|_{L_\iy(V)}=
\left\vert\Delta_N(ix_0)\right\vert$ for a point
$x_0\in V$. Given $a>0$, we define the set
$A(a):=\{x\in V: x=k/a,\,k\in\Z^m\}$.
Then $d(x_0,A(a))\le \sqrt{m}/a$. Therefore, for any $\vep>0$
there exist $a_0>0$ and $x_{0,a}=k_a/a\in A(a)$ such that for all
$a>a_0,\,\vert \Delta_N(ix_0)-\Delta_N(ix_{0,a}\vert <\vep$.
Hence
\beq\label{E2.12a}
\left\|\Delta_N(i\cdot)\right\|_{L_\iy(V)}
<\left\vert\Delta_N(ix_{0,a})\right\vert+\vep
\le a^{-N}\max_{k\in aV\cap\Z^m}\left\vert\Delta_N(ik)\right\vert+\vep
\le \left\|\Delta_N(i\cdot)\right\|_{L_\iy(V)}+\vep.
\eeq
Then \eqref{E1.7bn} follows from \eqref{E1.6n} and \eqref{E2.12a}.

Further, we assume that $a>1$. Let $Q_{M/a}$ be a cube with $M\in\N$,
satisfying the conditions $V\subset Q_{M/a}$ and
 $d\left(\partial Q_{M/a},V\right)>\sqrt{m}$.
It is easy to see that if $Q$ is a cube in $\R^m$ with an edge length
of $1/a$ such that $Q\cap V\ne\emptyset$, then $Q\subseteq Q_{M/a}$.
Next  denoting by  $\chi_{V}$ the characteristic function of the set $V$,
we see that the sum
\ba
\frac{1}{a^{2N+m}}\sum_{k\in aV\cap\Z^m}
\left\vert\Delta_N\left({ik}\right)\right\vert^2
=\frac{1}{a^{m}}\sum_{k\in aV\cap\Z^m}
\left\vert\Delta_N\left(\frac{ik}{a}\right)\right\vert^2
=\frac{1}{a^{m}}\sum_{k/a\in Q_{M/a}\cap\Z^m}
\left\vert\Delta_N\left(\frac{ik}{a}\right)\right\vert^2
\chi_V\left(\frac{k}{a}\right),
\ea
is the Riemann sum of the function
$\left\vert\Delta_N(i\cdot)\right\vert^2
\chi_V(\cdot):Q_{M/a}\to\R^1$, corresponding to the partition
\ba
\left\{\left\{x\in\R^m:\frac{k_j-1}{a}\le x_j\le \frac{k_j}{a},\,
1\le j\le m\right\}, \,k_j\in\Z^1,\,-M+1\le k_j\le M,\,1\le j\le m\right\}
\ea
of $Q_{M/a}$. Therefore, this sum converges to
$\int_V\left\vert\Delta_N(ix)\right\vert^2dx$ as $a\to \iy$.
This establishes
\eqref{E1.7an}.\hfill $\Box$

\section{Multivariate Levitan's Polynomials}\label{S2}
 \noindent
\setcounter{equation}{0}
\noindent
To prove Theorem \ref{T1.2}, we need a multivariate version of Levitan's
trigonometric polynomials introduced in the univariate case by
Levitan \cite{L1937} and H\"{o}rmander \cite{H1954, H1955}. Let us
define these polynomials and study their properties.

Let us set
\beq \label{E2.1n}
h(y):=\prod_{j=1}^m\frac{\sin y_j}{y_j},\qquad y\in\R^m.
\eeq
Then $h^2(\cdot/2)\in B_{Q_1}\cap L_1(\R^m)$ and
\beq \label{E2.2n}
\sum_{k\in\Z^m}h^2(x/2+k\pi)
=\sum_{k_1=-\iy}^{\iy}\left(\frac{\sin (x_1/2)}{x_1/2+k_1 \pi}\right)^2
\cdot\cdot\cdot
\sum_{k_m=-\iy}^{\iy}\left(\frac{\sin (x_m/2)}{x_m/2+k_m \pi}\right)^2=1,
\qquad x\in \R^m.
\eeq
Note that for $m=1$ relations \eqref{E2.2n} follows from a well-known expansion of
$\sin^{-2}(\tau/2),\,\tau\in\R^1$.

Let $f\in B_V\cap L_\iy(\R^m)$. Then for a fixed number $a>0$ the function
$g_a(u):=f(u)h^2(u/(2a))$
is an entire function and by Lemma \ref{L2.1}(b),
\ba
\vert g_a(x+iy)\vert
\le\|f\|_{L_\iy(\R^m)}\exp\left(\|y\|^*_V+a^{-1}\sum_{j=1}^m\vert y_j\vert\right)
\le\|f\|_{L_\iy(\R^m)}\exp((1+c/a)\|y\|^*_V),
\ea
where
\beq\label{E3.2a}
 c:=\sup_{y\in\R^m\setminus\{0\}}\frac{\sum_{j=1}^m\vert y_j\vert}{\|y\|^*_V}.
 \eeq
Therefore, $g_a\in B_{(1+c/a)V}\cap L_1(\R^m)$ and by
Nikolskii-type inequality \eqref{E1.9n}, $g_a\in L_2(\R^m)\cap B_{(1+c/a)V}$. Then by
Lemma \ref{L2.1}(c), $\mbox{supp}\,\widehat{g_a}\subseteq (1+c/a)V$.
Hence the Fourier series
\ba
(2\pi)^{-m/2}\sum_{k\in\Z^m}a^{-m}\widehat{g_a}(k/a)e^{ik\cdot x/a}
\ea
coincides with a trigonometric polynomial
\beq\label{E3.2b}
S_a(x)=S_a(f,x):=(2\pi)^{-m/2}\sum_{k\in (a+c)V\cap\Z^m}a^{-m}\widehat{g_a}(k/a)e^{ik\cdot x/a}
\eeq
of period $2\pi a$ in each variable with its spectrum in $(a+c)V$. This polynomial $S_a=S_a(f,\cdot)$
we call the \emph{multivariate Levitan's polynomial for} $f$.

For $m=1,\,V=[-\sa,\sa]$, and $a=n/\sa,\,n\in\N$, a trigonometric polynomial $S_a$ of
degree $n$ was introduced in \cite{L1937, H1954, H1955} (see also \cite[Sect. 85]{A1965}).
In the case when $V$ is a parallelohedron, multivariate Levitan's polynomials were
introduced in \cite{G1990}.

Let us discuss properties of $S_a$.

\begin{property} \label{P2.0}
If $f\in B_V\cap L_\iy(\R^m)$, then the following representation holds true:
\beq\label{E2.3n}
S_a(f,x)=\sum_{k\in\Z^m}f(x+2k\pi a)h^2(x/(2a)+k\pi),\qquad x\in\R^m,
\eeq
where $h$ is defined in \eqref{E2.1n}.
\end{property}
\proof
Setting $g(y):=f(ay) h^2(y/2)=g_a(ay)$, we see that $g\in L_1(\R^m)$
and the series
\beq\label{E3.3na}
\sum_{k\in\Z^m}\widehat{g}(k) e^{ik\cdot y}
=\sum_{k\in\Z^m}a^{-m}\widehat{g_a}(k/a)e^{ik\cdot y}
=(2\pi)^{m/2} S_a(ay)
\eeq
is a finite sum. Then using Lemma \ref{L2.2a},
we conclude that \eqref{E2.2a} and \eqref{E3.3na} imply \eqref{E2.3n} for $y=x/a$.
 \hfill $\Box$

\begin{property} \label{P2.1}
 If $f$ is a real-valued function from $B_V\cap L_\iy(\R^m)$, then $S_a(f,\cdot)$
 is a real-valued polynomial.
 \end{property}
 This property follows immediately from representation \eqref{E2.3n}.
 \begin{property} \label{P2.2}
 If $f\in B_V\cap L_p(\R^m),\,p\in(0,\iy]$, then
 \beq \label{E2.4n}
 \|S_a\|_{L_p(Q_{a\pi})}\le \|f\|_{L_p(\R^m)}.
 \eeq
 \end{property}
 \proof
 We first note that $f\in B_V\cap L_\iy(\R^m)$ by Nikolskii-type inequality \eqref{E1.9n},
 so by Property \ref{P2.0}, $S_a$ is represented by \eqref{E2.3n}. If $p\in(0,\iy)$,
 then
 \beq \label{E2.5n}
 \vert S_a(x)\vert^p\le \sum_{k=-\iy}^\iy \vert f(x+2k\pi a)\vert^p.
 \eeq
 Indeed, for $p=1$ \eqref{E2.5n} follows immediately from \eqref{E2.3n} since
 $\| h\|_{L_\iy(\R^m)}\le 1$.
 If $p\in (1,\iy)$, then by H\"{o}lder's inequality,  we obtain from relations
 \eqref{E2.3n} and \eqref{E2.2n} that
 \ba
  \vert S_a(x)\vert
  &\le& \left(\sum_{k\in\Z^m} \vert f(x+2k\pi a)\vert^p\right)^{1/p}
  \left(\sum_{k\in\Z^m}
  \left \vert h(x/(2a)+k\pi)\right\vert^{2p/(p-1)}\right)^{(p-1)/p}\\
  &\le& \left(\sum_{k\in\Z^m} \vert f(x+2k\pi a)\vert^p\right)^{1/p}.
  \ea
  If $p\in(0,1)$, then we have from \eqref{E2.3n}
  \ba
  \vert S_a(x)\vert^p\le\sum_{k\in\Z^m}
  \vert f(x+2k\pi a)\vert^p \left\vert h(x/(2a)+k\pi) \right\vert^{2p}
  \le \sum_{k\in\Z^m} \vert f(x+2k\pi a)\vert^p.
  \ea
  Thus \eqref{E2.5n} holds true for $p\in(0,\iy)$.  Next, integrating
  \eqref{E2.5n} over $Q_{a\pi}$
  and using Fatou's Lemma,
  we arrive at \eqref{E2.4n} for $p\in(0,\iy)$.
  Finally, for $p=\iy$ \eqref{E2.4n}
  follows immediately from \eqref{E2.3n} and \eqref{E2.2n}.\hfill $\Box$
 \vspace{.12in}\\

  Next, we discuss three approximation properties of
  multivariate Levitan's polynomials.
   \begin{property} \label{P2.3}
   If $f\in B_V\cap L_\iy(\R^m)$, then for $x\in\R^m$,
  \beq \label{E2.6n}
  \vert f(x)-S_a(x)\vert \le C_3 (\vert x\vert/a)^2\|f\|_{L_\iy(\R^m)},
  \eeq
  where $C_3\le 1/6$ is an absolute constant.
  \end{property}
  \proof
  To prove \eqref{E2.6n}, we need to evaluate certain integrals.
  First, we note that by Fubini's theorem,
  \ba
  \int_{Q_{1/2}}(x\cdot u/a)^{2l+1}du=0,\qquad a>0,\qquad l\in\Z^1_+,
  \ea
   because each monomial in the expansion of
  $(x\cdot u/a)^{2l+1}$ contains at least one factor of the form $u_j^{\al_j}$,
   where $\al_j>0$ is an odd number, $1\le j\le m$. Then
   \beq\label{E2.7nn}
   \int_{Q_{1/2}}\sin(x\cdot u/a)\,du
   =\sum_{l=0}^\iy\frac{(-1)^l}{(2l+1)!}
   \int_{Q_{1/2}}(x\cdot u/a)^{2l+1}du=0.
   \eeq
  Next, it follows from \eqref{E2.2n}, \eqref{E2.3n}, and
  \eqref{E2.7nn} that
  \bna
  &&\vert f(x)-S_a(x)\vert
  \le  \vert f(x)\vert \left(1-h^2(x/(2a))\right)
  +\sum_{k\ne 0}\vert f(x+2k\pi a)\vert
  h^2(x/(2a)+k\pi)\nonumber\\
  &&\le  2\left(1-\prod_{j=1}^m\left(\frac{\sin\left(\frac{x_j}{2a}\right)}{x_j/(2a)}
  \right)^2\right) \|f\|_{L_\iy(\R^m)}
  \le 4\left(1-\prod_{j=1}^m\frac{\sin\left(\frac{x_j}{2a}\right)}{x_j/(2a)}
  \right) \|f\|_{L_\iy(\R^m)}\nonumber\\
 &&=4\int_{Q_{1/2}}(1-\exp(i(x\cdot u/a)))\,du\, \|f\|_{L_\iy(\R^m)}
  = 4\int_{Q_{1/2}}(1-\cos(x\cdot u/a))\,du\, \|f\|_{L_\iy(\R^m)}\nonumber\\
  &&\le \frac{2}{a^2}\int_{Q_{1/2}} (x\cdot u)^2du\, \|f\|_{L_\iy(\R^m)}
  = \frac{\vert x\vert ^2}{6 a^2}\,\|f\|_{L_\iy(\R^m)}.
  \ena
 Thus \eqref{E2.6n} is established. \hfill $\Box$

   \begin{property} \label{P2.4}
   Let $f\in B_V\cap L_\iy(\R^m)$ and let $D_N:=\sum_{\vert \al\vert=N}b_\al D^\al$ be a
   linear differential operator with constant coefficients $b_\al\in\CC^1,\,\vert \al\vert=N$.
   Then for $x\in\R^m$ and $a\ge 1$,
   \beq \label{E2.7n}
  \vert D_N(f)(x)-D_N(S_a)(x)\vert
  \le C_4(\mathrm{diam}(V))^N\left(\frac{\vert x\vert^2}{a^2}+\frac{1}{a}\right)
  \|f\|_{L_\iy(\R^m)},
  \eeq
  where $C_4$ is independent of $x,\,a,\, f$, and $V$.
  \end{property}
  \proof
  For $N=0$ Property \ref{P2.4} follows from Property \ref{P2.3}, so we assume
  that $N\ge 1$. We also note that the uniform convergence of all series below
  follows easily from identity \eqref{E2.2n}, representation \eqref{E2.3n}, and
  Bernstein-type inequality \eqref{E1.11n} for $q=p$.
  Next, recalling notation $g_a(u)=f(u)h^2(u/(2a))$,  we see from \eqref{E2.3n} that
  \bna \label{E2.8n}
  &&\vert D_N(f)(x)-D_N(S_a)(x)\vert\nonumber\\
  &&\le \left\vert D_N(f)(x)-\sum_{k\in\Z^m}D_N(f)(x+2k\pi a)h^2(x/(2a)+k\pi)
  \right\vert
  \nonumber\\
  &&+\left\vert \sum_{k\in\Z^m}D_N(g_a)(x+2k\pi a)-D_N(f)(x+2k\pi a)h^2(x/(2a)+k\pi)
  \right\vert
  \nonumber\\
  &&=I_1(x)+I_2(x).
  \ena
  Since $D_N(f)\in B_V\cap L_\iy(\R^m)$ by Lemma \ref{L2.1}(d), Remark \ref{R2.1a},
  and Bernstein-type inequality
   \eqref{E1.11n} for $q=p=\iy$,
   we have from Property \ref{P2.3}
  \beq \label{E2.9n}
  I_1(x)=\vert D_N(f)(x)-S_a(D_N(f),x)\vert
  \le C_5(\mathrm{diam}(V))^N(\vert x\vert /a)^2 \|f\|_{L_\iy(\R^m)},
  \eeq
  where $C_5$ is independent of $x,\,a,\,f,$ and $V$.
  Further, using  Bernstein-type inequality \eqref{E1.11n} for $q=p$ again,
  we obtain by the multivariate Leibniz formula
  \bna \label{E2.10n}
  &&I_2(x)\nonumber\\
  &&\le \sum_{k\in\Z^m}\sum_{\vert\al\vert=N}\vert b_\al\vert
  \sum_{\be\le\al,\,\vert\be\vert\ge 1}\binom{\al}{\be}
  \left\vert D^{\al-\be}(f)(x+2k\pi a)\right\vert
  \left\vert D^{\be}(\g_a)(x+2k\pi a)\right\vert\nonumber\\
  &&\le C_6 (1+\mathrm{diam}(V))^N \|f\|_{L_\iy(\R^m)}
  \sum_{\vert\al\vert=N}\vert b_\al\vert
  \sum_{\be\le\al,\,\vert\be\vert\ge 1}
  \sum_{k\in\Z^m}
  \left\vert D^{\be}(\g_a)(x+2k\pi a)\right\vert
  \nonumber\\
  &&=C_6 (1+\mathrm{diam}(V))^N  \|f\|_{L_\iy(\R^m)} \nonumber\\
  &&\times\sum_{\vert\al\vert=N}\vert b_\al\vert
  \sum_{\be\le\al,\,\vert\be\vert\ge 1}
  (2a)^{-\vert\be\vert}
  \sum_{k\in\Z^m}
  \left\vert D^{\be}(\g_{1/2})(x/(2a)+k\pi)\right\vert,
  \ena
  where $\g_a(u):=h^2(u/(2a))$
  and $C_6$ is independent of $x,\,a,\,f,$ and $V$.
  It remains to estimate the series
  \beq \label{E2.11n}
  \sum_{k\in\Z^m}
  \left\vert D^{\be}(\g_{1/2})(y+2k\pi)\right\vert
  =\prod_{j=1}^m H_{\be_j}(y_j),\qquad y\in\R^m,
  \eeq
  where
  \beq \label{E2.12n}
  H_d(\tau)
  :=\sum_{l=-\iy}^{\iy}\left\vert\left(\left(\frac{\sin
  \tau}{\tau+l \pi}\right)^2\right)^{(d)}\right\vert
  \le C_7(d),\qquad \tau\in\R^1,
  \eeq
  by inequalities (2.12) and (2.14) in \cite{GT2017}. Therefore,
  by \eqref{E2.10n}, \eqref{E2.11n} and \eqref{E2.12n}, we obtain
  the estimate
  \beq \label{E2.13n}
  I_2(x)\le C_8 a^{-1} (1+\mathrm{diam}(V))^N \|f\|_{L_\iy(\R^m)},
  \eeq
  where $C_8$ is independent of $x,\,a,\,f,$ and $V$.
  Combining \eqref{E2.8n} with \eqref{E2.9n} and \eqref{E2.13n},
  we arrive at \eqref{E2.7n}.\hfill $\Box$

  \begin{property} \label{P2.5}
  Let $\g(a):[1,\iy)\to (0,1)$ be a function such that $\lim_{a\to\iy}\g(a)=0$.
   For any $V\subset \R^m,\,q\in(0,\iy]$, and $N\in\Z_+^1$ there exists a family of
   numbers $\{M_a\}_{a\in[1,\iy)}$, satisfying the following conditions:
   \begin{itemize}
   \item[(a)] $M_a\in(0,a\pi],\, a\in[1,\iy)$.
   \item[(b)] $\lim_{a\to\iy}M_a=\iy$.
   \item[(c)] For every family of functions $\{f_{(a)}\}_{a\in[1,\iy)}$
   such that
   $f_{(a)}\in B_{(1-\g(a))V}\cap L_\iy(\R^m),\,a\in[1,\iy)$, and
   $\sup_{a\in[1,\iy)}\left\|f_{(a)}\right\|_{L_\iy(\R^m)}<\iy,$
   the following relation holds true:
   \ba
   \lim_{a\to\iy}\|D_N\left(f_{(a)}\right)-D_N\left(S_a\left(f_{(a)},\cdot\right)\right)\|_{L_q(Q_{M_a})}=0.
   \ea
   \end{itemize}
   \end{property}
   \proof
   Let us set $M_a:=\min\{a^\de,a\pi \}$, where $\de\in(0,\vep_{q,N})$ and
  \ba
  \vep_{q,N}:=\left\{\begin{array}{ll}
  2q/(2q+m), &N=0,\\
  \min\{q/m,2q/(2q+m)\},&N\in\N,
  \end{array}\right.
  \ea
  for $q\in(0,\iy]$. Then conditions (a) and (b) are satisfied.
  Next by Property \ref{P2.3} for $N=0$ and by Property \ref{P2.4} for $N\in\N$,
  we obtain for $q\in(0,\iy]$
  \ba
  \|D_N\left(f_{(a)}\right)&-&D_N\left(S_a\left(f_{(a)},\cdot\right)\right)\|_{L_q(Q_{M_a})}\\
  &=&\left\{\begin{array}{ll}
  O([(1-\g(a))\mathrm{diam}(V)+1]^N M_a^{2+m/q}/a^2), &N=0,\\
  O([(1-\g(a))\mathrm{diam}(V)+1]^N M_a^{2+m/q}/a^2+M_a^{m/q}/a),&N\in\N,
  \end{array}\right.
  =o(1),
  \ea
 as $a\to\iy$. Hence condition (c) is satisfied as well.\hfill $\Box$

 \begin{remark}\label{R2.6}
  In case of $m=1$ and $\g(a)=0$, Property \ref{P2.0} was established in
 \cite{L1937, H1954, H1955}, while Properties \ref{P2.1} through \ref{P2.5}
  were proved in \cite{GT2017}.
 \end{remark}

\section{Proofs of Theorems}\label{S4}
 \noindent
\setcounter{equation}{0}
\emph{Proof of Theorem \ref{T1.2}.}
Let $f\in B_V\cap L_p(\R^m)$, where $p\in (0,\iy]$, and let $c$ be defined by
\eqref{E3.2a}. Setting $f_{(a)}(x):=f((1+c/a)^{-1}x),\,a\ge 1$, we see from Definition
\ref{D1.1} and Nikolskii-type inequality \eqref{E1.9n} that
$f_{(a)}\in B_{(1+c/a)^{-1}V}\cap L_\iy(\R^m)$. We can now
consider
the multivariate Levitan's polynomial $S_a\left(f_{(a)},\cdot\right)$ defined by
\eqref{E3.2b} and \eqref{E2.3n}. Since $S_a\left(f_{(a)},a\cdot\right)\in\TT_{aV}$
by \eqref{E3.2b},
we obtain for $0<p\le q\le\iy$
\beq \label{E4.1}
\|D_N(S_a(f_{(a)},\cdot)\|_{L_q(Q_{a\pi})}
\le P_{p,q,D_N,a,V} \|S_a(f_{(a)},\cdot)\|_{L_p(Q_{a\pi})}.
\eeq
In addition, we note that $D_N(f)\in L_q(\R^m)$ by Bernstein-Nikolskii type
inequality \eqref{E1.11n}. Therefore setting $\g(a):=c/(c+a)$, and using Properties
\ref{P2.5}, \ref{P2.2}, and inequalities \eqref{E1.1a} and \eqref{E4.1}, we obtain
 \ba
 \|D_N(f)\|_{L_q(\R^m)}
 &=&\lim_{a\to\iy}\|D_N(f)\|_{L_q\left((1+c/a)^{-1}Q_{M_a}\right)}\\
 &=&\lim_{a\to\iy}\|D_N\left(f_{(a)}\right)\|_{L_q\left(Q_{M_a}\right)}\\
 &\le & \left(\lim_{a\to\iy}\|D_N\left(f_{(a)}\right)
 -D_N\left(S_a\left(f_{(a)},\cdot\right)\right)\|_{L_q\left(Q_{M_a}\right)}^{\tilde{q}}\right.\\
 &+&\left.\liminf_{a\to\iy}\|D_N\left(S_a\left(f_{(a)},\cdot\right)\right)\|_{L_q\left(Q_{M_a}\right)}^{\tilde{q}}\right)^{1/\tilde{q}}\\
 &\le& \liminf_{a\to\iy}\|D_N\left(S_a\left(f_{(a)},\cdot\right)\right)\|_{L_q\left(Q_{a\pi}\right)}\\
 &\le& \liminf_{a\to\iy}\left(P_{p,q,D_N,a,V}\,
 \|S_a\left(f_{(a)},\cdot\right)\|_{L_p\left(Q_{a\pi}\right)}\right)\\
 &\le& \liminf_{a\to\iy}\left(P_{p,q,D_N,a,V} \,\|f_{(a)}\|_{L_p(\R^m)}\right)\\
 &=&\liminf_{a\to\iy}P_{p,q,D_N,a,V} \,\|f\|_{L_p(\R^m)}.
 \ea
 Thus \eqref{E1.12n} is established.\hfill $\Box$ \vspace{.1in}\\
\emph{Proof of Theorem \ref{T1.3}.}
We first note that there exists $a_0\in(0,\iy)$ such that for
$p\in(0,\iy],\,
N\in\Z^m,$ and $V\subset\R^m$,
the following crude estimates of $P_{p,\iy,D_N,a,V}$ hold true:
\beq \label{E4.1a}
C_9(p,D_N,V)
\le \inf_{a\in[a_0,\iy)}P_{p,\iy,D_N,a,V}
\le\sup_{a\in[a_0,\iy)}P_{p,\iy,D_N,a,V}\le C_{10}(p,D_N,V).
\eeq
The left and right inequalities in \eqref{E4.1a} follow from \eqref{E1.12n}
and \eqref{E1.10n}, respectively.

We will prove the
 theorem by constructing a nontrivial function $f_0\in B_V\cap L_p(\R^m),\,
 p\in(0,\iy]$, such that
 \beq \label{E4.2}
E_{p,\iy,D_N,V} \ge \|D_N\left(f_0\right)\|_{L_\iy(\R^m)}/
\|f_0\|_{L_p(\R^m)}
 \ge \limsup_{a\to\iy}P_{p,\iy,D_N,a,V}.
 \eeq
 Then combining \eqref{E4.2} with \eqref{E1.12n} for $q=\iy$, we see that
$\lim_{a\to\iy}P_{p,\iy,D_N,a,V}$ exists and
\beq \label{E4.3}
E_{p,\iy,D_N,V}=\lim_{a\to\iy}P_{p,\iy,D_N,a,V}.
\eeq
This proves \eqref{E1.13n}.
In addition, $f_0$ is an extremal function in \eqref{E4.3};
that is, \eqref{E1.14n} is valid.

It remains to construct a function $f_0$, satisfying \eqref{E4.2}.
Let $T_a\in\TT_{aV}$ be a polynomial, satisfying the equality
\beq \label{E4.4}
P_{p,\iy,D_N,a,V}=a^{-N-m/p}\|D_N\left(T_a\right)\|_{L_\iy\left(Q_\pi\right)}
/\|T_a\|_{L_p\left(Q_\pi\right)}.
\eeq
The existence of an extremal polynomial $T_a$ in \eqref{E4.4}
can be proved by the standard compactness argument. Indeed,
given $d\in\N$, let $T_{a,d}\in\TT_{aV}$ satisfy the following
relations:
\beq\label{E4.4a}
\|T_{a,d}\|_{L_\iy\left(Q_\pi\right)}=1,\qquad
\|D_N\left(T_{a,d}\right)\|_{L_\iy\left(Q_\pi\right)}
=\left\vert D_N\left(T_{a,d}\right)(0)
\right\vert,
\eeq
and
\ba
a^{N+m/p}P_{p,\iy,D_N,a,V}
<\|D_N\left(T_{a,d}\right)\|_{L_\iy\left(Q_\pi\right)}
/\|T_{a,d}\|_{L_p\left(Q_\pi\right)} +1/d.
\ea
Then there exists a nontrivial polynomial $T_a\in\TT_{aV}$ and a sequence
$\{d_k\}_{k=1}^\iy\subseteq\N$ such that for any $\al\in\Z^m_+$,
$\lim_{k\to\iy}D^\al T_{a,d_k}(x)=D^\al T_a(x)$
uniformly on $Q_\pi$. Thus \eqref{E4.4} holds true.

Next setting $U_a(x):=T_a(x/a)$, we see that
$U_a\in B_V\cap L_\iy(\R^m)$. In addition, it follows from \eqref{E4.4}
that
\beq \label{E4.5}
P_{p,\iy,D_N,a,V}=\|D_N\left(U_a\right)\|_{L_\iy \left(Q_{a\pi}\right)}
/\|U_a\|_{L_p \left(Q_{a\pi}\right)}.
\eeq
Moreover, we can assume that
\beq \label{E4.6}
\|D_N\left(U_a\right)\|_{L_\iy \left(Q_{a\pi}\right)}=\left\vert D_N\left(U_a\right)(0)
\right\vert=1.
\eeq
Note that normalization \eqref{E4.6} is different compared with \eqref{E4.4a}.
Normalization \eqref{E4.4a} was used only for the proof of the existence
of extremal polynomials $T_a$ in \eqref{E4.4}.

Then we obtain from \eqref{E4.5}, \eqref{E4.6}, and \eqref{E4.1a}
\bna \label{E4.7}
&&\sup_{a\in[a_0,\iy)} \|U_a\|_{L_\iy(\R^m)}
\le \sup_{a\in[a_0,\iy)}\left(P_{p,\iy,D_0,a,V}\|U_a\|_{L_p\left(Q_{a\pi}\right)}
\right)\nonumber\\
&=&\sup_{a\in[a_0,\iy)}\left(P_{p,\iy,D_0,a,V}
/P_{p,\iy,D_N,a,V}\right)
\le C_{10}(p,D_0,V)/C_9(p,D_N,V).
\ena
Let $\{a_s\}_{s=1}^\iy\subset [a_0,\iy)$ be a sequence such that
\beq \label{E4.8}
\limsup_{a\to\iy}P_{p,\iy,D_N,a,V}=\lim_{s\to\iy}P_{p,\iy,D_N,a_s,V}.
\eeq
Using now the compactness theorem of Lemma \ref{L2.2} for the sequence
$\{U_{a_s}\}_{s=1}^\iy$ of functions from $B_V\cap L_\iy(\R^m)$
uniformly bounded by \eqref{E4.7}, we see that there exist
 a function $f_0\in B_V\cap L_\iy(\R^m)$ and a subsequence
 $\{a_{s_r}\}_{r=1}^\iy$ such that
 \beq \label{E4.9}
 \lim_{r\to\iy}D^\al U_{a_{s_r}}(x)=D^\al f_0(x),\qquad \al\in\Z^m_+,
 \eeq
 uniformly on any cube $Q_A\subset\R^m,\,A>0$. Moreover, by
  \eqref{E4.6} and \eqref{E4.9},
  \beq \label{E4.10}
  \|D_N\left(f_{0}\right)\|_{L_\iy(\R^m)}=\left\vert D_N\left(f_{0}\right)(0)
\right\vert=1.
\eeq
In addition, using \eqref{E1.1a}, \eqref{E4.9}, \eqref{E4.5}, and \eqref{E4.6},
we obtain for any cube $Q_A\subset\R^m,\,A>0$,
\bna \label{E4.11}
&&\|f_0\|_{L_p\left(Q_A\right)}^{\tilde{p}}
\le \lim_{r\to\iy}
\left(\|f_0-U_{a_{s_r}}\|_{L_p\left(Q_A\right)}^{\tilde{p}}
+\|U_{a_{s_r}}\|_{L_p\left(Q_A\right)}^{\tilde{p}}\right)\nonumber\\
&&\le \lim_{r\to\iy}\|U_{a_{s_r}}\|^{\tilde{p}}
_{L_p\left(Q_{a_{s_r}\pi}\right)}
=1/
\lim_{r\to\iy}P_{p,\iy,D_N,a_{s_r},V}^{\tilde{p}}.
\ena
Next using \eqref{E4.11} and \eqref{E4.1a}, we see that
\beq \label{E4.12}
\|f_0\|_{L_p\left(Q_A\right)}\le 1/C_9(p,D_N,V).
\eeq
Therefore, $f_0$ is a nontrivial function from $B_V\cap L_p(\R^m)$,
by \eqref{E4.10} and \eqref{E4.12}. Thus for any cube
$Q_A\subset\R^m,\,A>0$, we obtain from \eqref{E4.8}, \eqref{E4.5}, \eqref{E4.9},
 and \eqref{E4.10}
\bna \label{E4.13}
&&\limsup_{a\to\iy}P_{p,\iy,D_N,a,V}
=\lim_{r\to\iy}\left(\|U_{a_{s_r}}\|_{L_p\left(Q_{a_{s_r}\pi}\right)}
\right)^{-1}\nonumber\\
&&\le \lim_{r\to\iy}\left(\|U_{a_{s_r}}\|_{L_p\left(Q_{A}\right)}
\right)^{-1}
=\|D_N\left(f_0\right)\|_{L_\iy(\R^m)}/
\|f_0\|_{L_p\left(Q_{A}\right)}.
\ena
Finally, letting $A\to \iy$ in \eqref{E4.13}, we arrive at \eqref{E4.2}.
\hfill$\Box$ \vspace{.1in}\\

\end{document}